\documentclass{myelsart}

\usepackage{graphicx}
\usepackage{amssymb}

\usepackage[english,francais]{babel}

\newtheorem{theorem}{Theorem}[section]
\newtheorem{lemma}[theorem]{Lemma}
\newtheorem{e-proposition}[theorem]{Proposition}

\newtheorem{e-definition}[theorem]{Definition\rm}

\newtheorem{theoreme}{Th\'eor\`eme}[section]

\newtheorem{proposition}[theoreme]{Proposition}

\renewcommand{\theequation}{\arabic{equation}}
\def\og{\leavevmode\raise.3ex\hbox{$\scriptscriptstyle\langle\!\langle$~}}
\def\fg{\leavevmode\raise.3ex\hbox{~$\!\scriptscriptstyle\,\rangle\!\rangle$}}

\journal{the Acad\'emie des sciences}
\begin{document}
% place in the next line the header (rubrique) chosen for your article,
% if you know it (you can also have 2, format : Header1/Header2
\centerline{}
\begin{frontmatter}


% \author{Name\thanksref{label2}}
\selectlanguage{english}
\title{The Boltzmann-Grad limit of the periodic Lorentz gas in two space dimensions}

% The [label1] can be suppressed if there is only one address for all authors

\selectlanguage{english}
\author[E. C.]{Emanuele Caglioti},
%\ead{caglioti@mat.uniroma1.it}
\author[F. G.]{Fran\c cois Golse}
%\ead{golse@math.polytechnique.fr}

\address[E. C.]{Universit\`a di Roma ``La Sapienza", Dipartimento di Matematica 
``Guido Castelnuovo", P.le Aldo Moro 2, 00185 Rome}
\address[F. G.]{Ecole polytechnique, Centre de Math\'ematiques Laurent Schwartz, 91128 Palaiseau Cedex}

% If you know the dates of reception, and acceptation you can put them now;
%  idem the name of the person presenting the Note

\medskip
\begin{center}
%{\small Received *****; accepted after revision +++++\\
%Presented by £££££}
\end{center}

\begin{abstract}
\selectlanguage{english}
% Text of abstract in English
The periodic Lorentz gas is the dynamical system corresponding to the free motion of a 
point particle in a periodic system of fixed spherical obstacles of radius $r$ centered at
the integer points, assuming all collisions of the particle with the obstacles to be elastic.
In this Note, we study this motion on time intervals of order $1/r$ as 
$r\to 0^+$.
%{\it To cite this article: E. Caglioti, F. Golse, C. R. Acad. Sci. Paris, Ser. I ??? (2007).}

\vskip 0.5\baselineskip

\selectlanguage{francais}
% Text of abstract in French
\noindent{\bf R\'esum\'e} \vskip 0.5\baselineskip \noindent
{\bf La limite de Boltzmann-Grad du gaz de Lorentz p\'eriodique en dimension deux d'espace. }
Le gaz de Lorentz p\'eriodique est le syst\`eme dynamique correspondant au mouvement 
libre dans le plan d'une particule ponctuelle rebondissant de mani\`ere \'elastique sur un syst\`eme de disques de rayon $r$ centr\'es aux points de coordonn\'ees enti\`eres. On
\'etudie ce mouvement pour $r\to 0^+$ sur des temps de l'ordre de $1/r$.
%{\it Pour citer cet article~: E. Caglioti, F. Golse, C. R. Acad. Sci. Paris, Ser. I ??? (2007).}

\end{abstract}
\end{frontmatter}

% now the Version française abrégée, if it exists
\selectlanguage{francais}
\section*{Version fran\c{c}aise abr\'eg\'ee}
% Text of your Version française abrégée here.
% Note you do not need to repeat here equations that you use in the
% main text - for example 'voir (3)' is quite acceptable.

On appelle gaz de Lorentz le syst\`eme dynamique correspondant au mouvement 
libre d'une particule ponctuelle dans un syst\`eme d'obstacles circulaires de rayon
$r$ centr\'es aux sommets d'un r\'eseau de $\mathbf{R}^2$, supposant que 
les collisions entre la particule et les obstacles sont parfaitement \'elastiques. Les 
trajectoires de la particule sont alors donn\'ees par les formules (\ref{Traj_r}).
La limite de Boltzmann-Grad pour le gaz de Lorentz consiste \`a supposer que le
rayon des obstacles $r\to 0^+$, et \`a observer la dynamique de la particule sur des 
plages de temps longues, de l'ordre de $1/r$ --- voir (\ref{Def-f_r}) pour la loi 
d'\'echelle de Boltzmann-Grad en dimension $2$.

Or les trajectoires de la particule s'expriment en fonction de l'application de transfert 
d'obstacle \`a obstacle $T_r$ d\'efinie par (\ref{TransitMap}) --- o\`u la notation $Y$ 
d\'esigne la transformation inverse de (\ref{Def-Y}) --- application qui associe, \`a tout 
param\`etre d'impact $h'\in[-1,1]$ correspondant \`a une particule quittant la surface 
d'un obstacle dans la direction $\omega\in\mathbf{S}^1$, le param\`etre d'impact 
$h$ \`a la collision suivante, ainsi que le temps $s$ s'\'ecoulant jusqu'\`a cette collision. 
(Pour une d\'efinition de la notion de param\`etre d'impact, voir (\ref{Tau-h}).)

On se ram\`ene donc \`a \'etudier le comportement limite de l'application de transfert 
$T_r$ pour $r\to 0^+$. 

\begin{proposition}
Lorsque $0<\omega_2<\omega_1$ et 
$\alpha=\frac{\omega_2}{\omega_1}\notin\mathbf{Q}$, l'application de transfert $T_r$ 
est approch\'ee \`a $O(r^2)$ pr\`es par l'application $\mathbf{T}_{A,B,Q,N}$ d\'efinie 
\`a la formule (\ref{TransitLimit}). Pour $\omega\in\mathbf{S}^1$ quelconque, on se 
ram\`ene au cas ci-dessus par la sym\'etrie (\ref{Sym}).  
\end{proposition}

Les param\`etres $A,B,Q,N\,\hbox{Êmod. }2$ intervenant dans l'application de transfert 
asymptotique sont d\'efinis \`a partir du d\'eveloppement en fraction continue 
(\ref{FracCont}) de $\alpha$ par les formules (\ref{Fla-Nk}) et (\ref{Fla-ABQ}). 

On voit sur ces formules que les param\`etres $A,B,Q,N\,\hbox{Êmod. }2$ sont des 
fonctions tr\`es fortement oscillantes des variables $\omega$ et $r$. Il est donc naturel 
de chercher le comportement limite de l'application de transfert $T_r$ dans une topologie 
faible vis \`a vis de la d\'ependance en la direction $\omega$. On montre ainsi que, pour 
tout $h'\in[-1,1]$, la famille d'applications $\omega\mapsto T_r(h',\omega)$ converge au 
sens des mesures de Young (voir par exemple \cite{Tartar} p. 146--154 pour une d\'efinition 
de cette notion) lorsque $r\to 0^+$ vers une mesure de probabilit\'e $P(s,h|h')dsdh$
ind\'ependante de $\omega$:

\begin{theoreme}
Pour tout $\Phi\in C_c(\mathbf{R}_+^*\times]-1,1[)$ et tout $h'\in]-1,1[$, la limite
(\ref{TransitProba}) a lieu dans $L^\infty(\mathbf{S}^1)$ faible-* lorsque $r\to 0^+$, o\`u la 
mesure de probabilit\'e $P(s,h|h')dsdh$ est l'image de la probabilit\'e $\mu$ d\'efinie dans
(\ref{PartiProba}) par l'application $(A,B,Q,N)\mapsto\mathbf{T}_{A,B,Q,N}(h')$ de la 
formule (\ref{TransitLimit}). De plus, cette densit\'e de probabilit\'e de transition $P(s,h|h')$
v\'erifie les propri\'et\'es (\ref{PropP}).
\end{theoreme}

Le th\'eor\`eme ci-dessus est le r\'esultat principal de cette Note: il montre que, dans la
limite de Boltzmann-Grad, le transfert d'obstacle \`a obstacle est d\'ecrit de mani\`ere
naturelle par une densit\'e de probabilit\'e de transition $P(s,h|h')$, o\`u $s$ est le laps 
de temps entre deux collisions successives avec les obstacles (dans l'\'echelle de temps
de la limite de Boltzmann-Grad), $h$ le param\`etre d'impact lors de la collision future
et $h'$ celui correspondant \`a la collision pass\'ee.

Le fait que la probabilit\'e de transition $P(s,h|h')$ soit ind\'ependante de la direction
sugg\`ere l'hypoth\`ese d'ind\'ependance (H) des quantit\'es $A,B,Q,N\,\hbox{Êmod. 2}$ 
correspondant \`a des collisions successives. 

\begin{theoreme}
Sous l'hypoth\`ese (H), 
pour toute densit\'e de probabilit\'e $f^{in}\in C_c(\mathbf{R}^2\times\mathbf{S}^1)$, la 
fonction de distribution $f_r\equiv f_r(t,x,\omega)$ de la th\'eorie cin\'etique, d\'efinie par 
(\ref{Def-f_r}) converge dans 
$L^\infty(\mathbf{R}_+\times\mathbf{R}^2\times\mathbf{S}^1)$ vers la limite (\ref{LimBG})
lorsque $r\to 0^+$, o\`u $F$ est la solution du probl\`eme de Cauchy (\ref{EqLim}) pos\'e
dans l'espace des phases \'etendu 
$(x,\omega,s,h)\in\mathbf{R}^2\times\mathbf{S^1}\times\mathbf{R}^*_+\times]-1,1[$.
\end{theoreme}

Dans le cas d'obstacles al\'eatoires ind\'ependants et poissonniens, Gallavotti a montr\'e 
que la limite de Boltzmann-Grad du gaz de Lorentz ob\'eit \`a l'\'equation cin\'etique 
de Lorentz (\ref{LorentzKinEq}). Le cas p\'eriodique est absolument diff\'erent: en se basant 
sur des estimations (cf. \cite{BourgainGolseWennberg} et \cite{GolseWennberg}) du temps 
de sortie du domaine $Z_r$ d\'efini dans (\ref{BillTable}), on d\'emontre que la limite de 
Boltzmann-Grad du gaz de Lorentz p\'eriodique ne peut pas \^etre d\'ecrite par l'\'equation 
de Lorentz (\ref{LorentzKinEq}) sur l'espace des phases $\mathbf{R}^2\times\mathbf{S}^1$
classique de la th\'eorie cin\'etique: voir \cite{Golse}. Si l'hypoth\`ese (H) ci-dessous \'etait
v\'erifi\'ee, le mod\`ele cin\'etique (\ref{LimBG}) dans l'espace des phases \'etendu fournirait
donc l'\'equation devant remplacer l'\'equation cin\'etique classique de Lorentz 
(\ref{LorentzKinEq}) dans le cas p\'eriodique.

\selectlanguage{english}
% main text
\section{The Lorentz gas}
\label{S-LorGas}
%%%%%%%%%%%%%%

The Lorentz gas is the dynamical system corresponding to the free motion of a single 
point particle in a periodic system of fixed spherical obstacles, assuming that collisions
between the particle and any of the obstacles are elastic. Henceforth, we assume that
the space dimension is $2$ and that the obstacles are disks of radius $r$ centered at
each point of $\mathbf{Z}^2$. Hence the domain left free for particle motion is
\begin{equation}
\label{BillTable}
Z_r=\{x\in\mathbf{R}^2\,|\,\hbox{dist}(x,\mathbf{Z}^2)>r\}\,,
	\qquad\hbox{ where it is assumed that $0<r<\frac12$.}
\end{equation}
Assuming that the particle moves at speed $1$, its trajectory starting from $x\in Z_r$ 
with velocity $\omega\in\mathbf{S}^1$ at time $t=0$ is 
$t\mapsto(X_r,\Omega_r)(t;x,\omega)\in\mathbf{R}^2\times\mathbf{S}^1$ given by
\begin{equation}
\label{Traj_r}
\begin{array}{lll}
\dot{X_r}(t)=\Omega_r(t)&\hbox{ and }\dot\Omega_r(t)=0
	&\hbox{ whenever $X_r(t)\in Z_r$,}
\\
X_r(t+0)=X_r(t-0)&\hbox{ and }\Omega_r(t+0)=\mathcal{R}[X_r(t)]\Omega_r(t-0)
	&\hbox{ whenever $X_r(t-0)\in\partial Z_r$,}
\end{array}\end{equation}
denoting $\dot{}=\frac{d}{dt}$ and $\mathcal{R}[X_r(t)]$ the specular reflection on 
$\partial Z_r$ at the point
$X_r(t)=X_r(t\pm 0)$. Assume that the initial position $x$ and direction $\omega$ of 
the particle are distributed in $Z_r\times\mathbf{S}^1$ with some probability density 
$f^{in}\equiv f^{in}(x,\omega)$, and define
\begin{equation}
\label{Def-f_r}
f_r(t,x,\omega):=f^{in}(rX_r(-t/r;x,\omega),\Omega_r(-t/r;x,\omega))
	\quad\hbox{ whenever $x\in Z_r$.}
\end{equation}
We are concerned with the limit of $f_r$ as $r\to 0^+$ in some appropriate sense to
be explained below. In the 2-dimensional setting considered here, this is precisely 
the Boltzmann-Grad limit.

In the case of a random (Poisson), instead of periodic, configuration of obstacles, 
Gallavotti \cite{Gallavotti} proved that the 
expectation of $f_r$ converges to the solution of the Lorentz kinetic equation for
$(x,\omega)\in\mathbf{R}^2\times\mathbf{S}^1$:
\begin{equation}
\label{LorentzKinEq}
(\partial_t+\omega\cdot\nabla_x)f(t,x,\omega)
	=\int_{\mathbf{S}^1}(f(t,x,\omega-2(\omega\cdot n)n)
		-f(t,x,\omega))(\omega\cdot n)_+dn\,,\qquad f\Big|_{t=0}=f^{in}\,.
\end{equation}
In the case of a periodic distribution of obstacles, the Boltzmann-Grad limit of the 
Lorentz gas cannot be described by a transport equation as above: see \cite{Golse} 
for a complete proof, based on estimates on the free path length to be found in 
\cite{BourgainGolseWennberg} and \cite{GolseWennberg}. This limit involves
instead a linear Boltzmann equation on an extended phase space with two new
variables taking into account correlations between consecutive collisions with
the obstacles that are an effect of periodicity: see Theorem \ref{T-LimEq}.

\section{The transfer map}
\label{S-TRANSIT}
%%%%%%%%%%%%%%%%

Denote by $n_x$ the inward unit normal to $Z_r$ at the point $x\in\partial Z_r$, 
consider 
\begin{equation}
\Gamma_r^\pm
	=\{(x,\omega)\in\partial Z_r\times\mathbf{S}^1\,|\,\pm\omega\cdot n_x>0\}\,, 
\end{equation}
and let $\Gamma^\pm_r/\mathbf{Z}^2$ be the quotient of $\Gamma^\pm_r$ 
under the action of $\mathbf{Z}^2$ by translation on the $x$ variable. For 
$(x,\omega)\in\Gamma^+_r$, let $\tau_r(x,\omega)$ be the exit time 
from $x$ in the direction $\omega$ and $h_r(x,\omega)$ be the impact 
parameter:
\begin{equation}
\label{Tau-h}
\tau_r(x,\omega)=\inf\{t>0\,|\,x+t\omega\in\partial Z_r\}\,,
	\quad\hbox{ and }h_r(x,\omega)=\sin(\omega,n_x)\,.
\end{equation}
Obviously, the map 
\begin{equation}
\label{Def-Y}
\Gamma^+_r/\mathbf{Z}^2\ni(x,\omega)
	\mapsto(h_r(x,\omega),\omega)\in]-1,1[\times\mathbf{S}^1
\end{equation} 
coordinatizes $\Gamma^+_r/\mathbf{Z}^2$, and we henceforth denote $Y_r$
its inverse. 

For each $r\in]0,\frac12[$, consider now the transfer map 
$T_r:\,]-1,1[\times\mathbf{S}^1\to\mathbf{R}_+^*\times]-1,1[$ defined by
\begin{equation}
\label{TransitMap}
T_r(h',\omega)=
(r\tau_r(Y_r(h',\omega)),h_r(X_r(\tau_r(Y_r(h',\omega));Y_r(h',\omega)),
	\Omega_r(\tau_r(Y_r(h',\omega));Y_r(h',\omega))))\,.
\end{equation}

For a particle leaving the surface of an obstacle in the direction $\omega$ with impact 
parameter $h'$, the transition map $T_r(h',\omega)=(s,h)$ gives the (rescaled) distance
$s$ to the next collision, and the corresponding impact parameter $h$. Obviously, each
trajectory (\ref{Traj_r}) of the particle can be expressed in terms of the transfer map $T_r$ 
and iterates thereof. The Boltzmann-Grad limit of the periodic Lorentz gas is therefore
reduced to computing the limiting behavior of $T_r$ as $r\to 0^+$, and this is our main
purpose in this Note.  

\smallskip

\smallskip
\begin{figure}
\centering

\includegraphics[width=6.0cm]{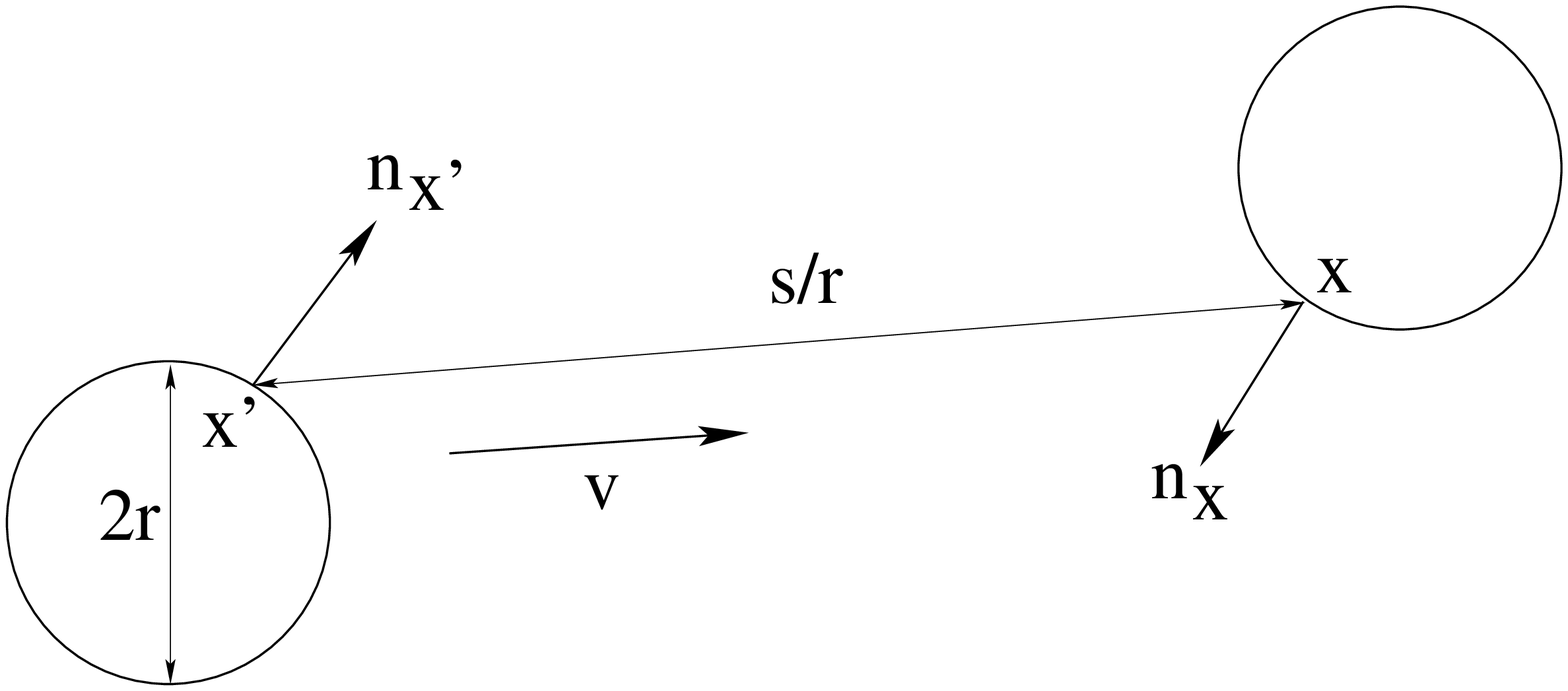}\hskip 2cm
\includegraphics[width=6.0cm]{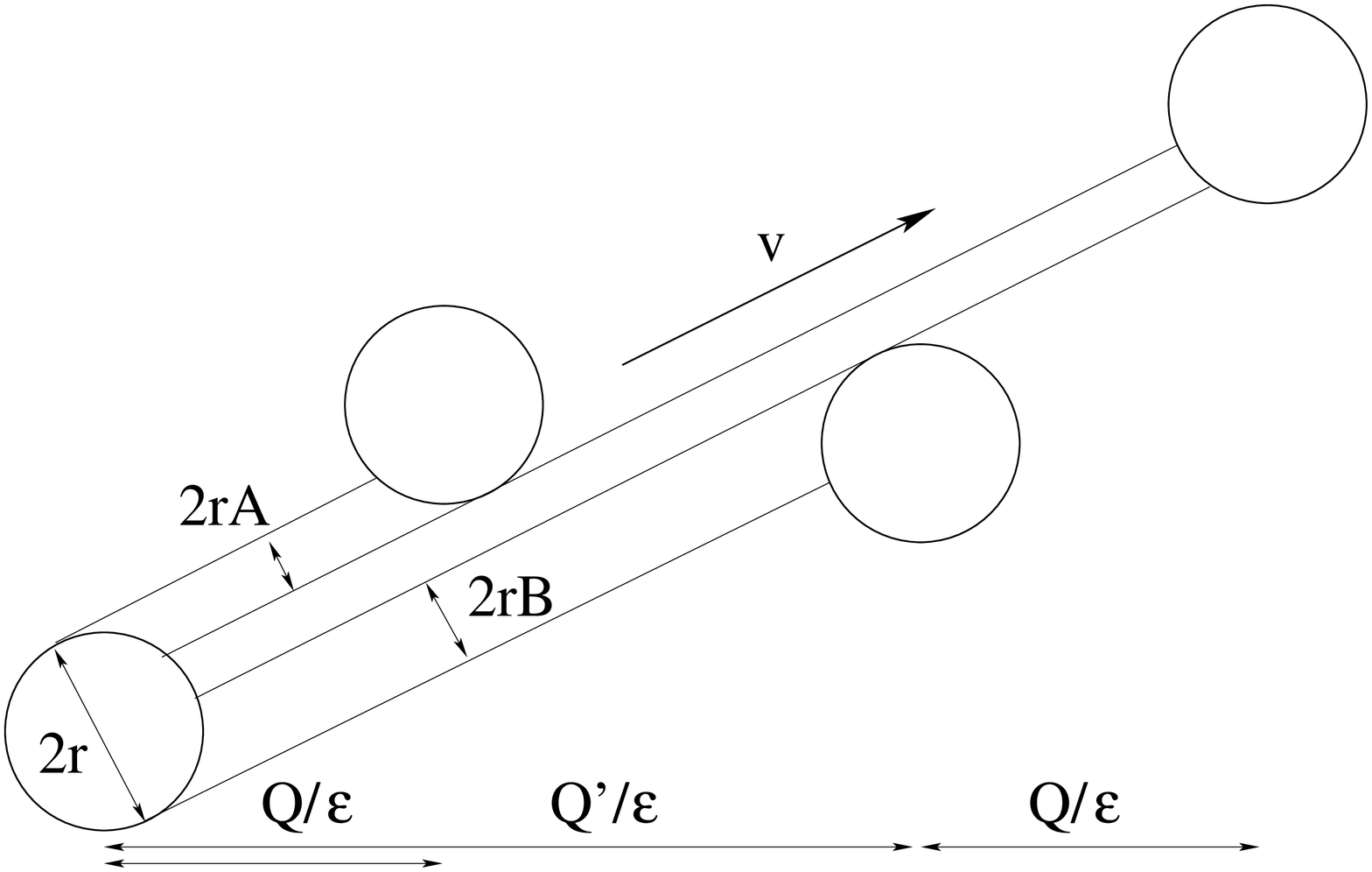}

\caption{Left: the transfer map $(s,h)=T_r(h',v)$, with $h'=\sin(n_{x'},v)$ and $h=\sin(n_x,v)$. Right:
Particles leaving the surface of one obstacle will next collide with one of generically 
three obstacles. The figure explains the geometrical meaning of $A,B,Q$.}

\end{figure}

\smallskip

We first need some pieces of notation. Assume $\omega=(\omega_1,\omega_2)$ with 
$0<\omega_2<\omega_1$, and $\alpha=\omega_2/\omega_1\in]0,1[\setminus\mathbf{Q}$. 
Consider the continued fraction expansion of $\alpha$:
\begin{equation}
\label{FracCont}
\alpha=[0;a_0,a_1,a_2,\ldots]
	=\frac1{\displaystyle a_0+\frac1{\displaystyle a_1+\ldots}}\,.
\end{equation}
Define the sequences of convergents $(p_n,q_n)_{n\ge 0}$ and errors 
$(d_n)_{n\ge 0}$ by the recursion formulas
\begin{equation}
\begin{array}{ll}
p_{n+1}=a_np_n+p_{n-1}\,,\quad &p_0=1\,,\,\,p_1=0\,,
	\qquad d_n=(-1)^{n-1}(q_n\alpha-p_n)\,,
\\
q_{n+1}=a_nq_n+q_{n-1}\quad &q_0=0\,,\,\,q_1=1\,,
\end{array}
\end{equation}
and let
\begin{equation}
\label{Fla-Nk}
N(\alpha,r)=\inf\{n\ge 0\,|\,d_n\le 2r\sqrt{1+\alpha^2}\}\,,\quad
\hbox{  and }k(\alpha,r)=
-\left[\frac{2r\sqrt{1+\alpha^2}-d_{N(\alpha,r)-1}}{d_{N(\alpha,r)}}\right]\,.
\end{equation}

\begin{e-proposition} 
\label{P-TransitMap}
For each $\omega=(\cos\theta,\sin\theta)$ with $0<\theta<\frac{\pi}4$, set
$\alpha=\tan\theta$ and $\epsilon=2r\sqrt{1+\alpha^2}$,  and 
\begin{equation}
\label{Fla-ABQ}
A(\alpha,r)=1-\frac{d_{N(\alpha,r)}}{\epsilon}\,,\quad
B(\alpha,r)=
1-\frac{d_{N(\alpha,r)-1}-k(\alpha,r)d_{N(\alpha,r)}}{\epsilon}\,,\quad
Q(\alpha,r)=\epsilon q_{N(\alpha,r)}\,.
\end{equation}
In the limit $r\to 0^+$, the transition map $T_r$ defined in (\ref{TransitMap}) is 
explicit in terms of $A,B,Q,N$ up to $O(r^2)$:
\begin{equation}
T_r(h',\omega)=
\mathbf{T}_{A(\alpha,r),B(\alpha,r),Q(\alpha,r),N(\alpha,r)}(h')+(O(r^2),0)
\hbox{ for each $h'\in]-1,1[$.}
\end{equation}
In the formula above
\begin{equation}
\label{TransitLimit}
\begin{array}{ll}
\mathbf{T}_{A,B,Q,N}(h')=(Q,h'-2(-1)^N(1-A))&\quad\hbox{ if }(-1)^Nh'\in]1-2A,1]\,,
\\
\mathbf{T}_{A,B,Q,N}(h')=\left(Q',h'+2(-1)^N(1-B)\right)
	&\quad\hbox{ if }(-1)^Nh'\in[-1,-1+2B[\,,
\\
\mathbf{T}_{A,B,Q,N}(h')=\left(Q'+Q,h'+2(-1)^N(A-B)\right)
	&\quad\hbox{ if }(-1)^Nh'\in[-1+2B,1-2A]\,,
\end{array}
\end{equation}
for each $(A,B,Q,N)\in K:=]0,1[^3\times\mathbf{Z}/2\mathbf{Z}$, with the notation
$Q'=\frac{1-Q(1-B)}{1-A}$.
\end{e-proposition}

The proof uses the 3-term partition of the 2-torus defined in section 2 of
\cite{CagliotiGolse}, following the work of \cite{BlankKrikorian}. 

For $\omega=(\cos\theta,\sin\theta)$ with arbitrary $\theta\in\mathbf{R}$, the map 
$h'\mapsto T_r(h',\omega)$ is computed using Proposition \ref{P-TransitMap} in 
the following manner. 
Set $\tilde\theta=\theta-m\frac{\pi}2$ with $m=[\frac2{\pi}(\theta+\frac{\pi}4)]$ and 
let $\tilde\omega=(\cos\tilde\theta,\sin\tilde\theta)$. Then
\begin{equation}
\label{Sym}
T_r(h',\omega)=(s,h)\,,\quad\hbox{ where }
(s,\hbox{sign}(\tan\tilde\theta)h)=T_r(\hbox{sign}(\tan\tilde\theta)h',\tilde\omega)\,.
\end{equation}

\section{The Boltzmann-Grad limit of the transfer map $T_r$}
\label{S-BGTM}
%%%%%%%%%%%%%%%%%%%%%%%%%%%%%%%%%%

The formulas (\ref{Fla-Nk}) and (\ref{Fla-ABQ}) defining $A,B,Q,N\,\hbox{Êmod. 2}$
show that these quantities are strongly oscillating functions of the variables $\omega$
and $r$. In view of Proposition \ref{P-TransitMap}, one therefore expects the transfer 
map $T_r$ to have a limit as $r\to 0^+$ only in the weakest imaginable sense, i.e. in 
the sense of Young measures --- see \cite{Tartar}, pp. 146--154 for a definition of this 
notion of convergence. 

The main result in the present Note is the theorem below. It says that, for each 
$h'\in[-1,1]$, the family of maps $\omega\mapsto T_r(h',\omega)$ converges 
as $r\to 0^+$ and in the sense of Young measures to some probability measure 
$P(s,h|h')dsdh$ that is moreover independent of $\omega$.

\begin{theorem}
\label{T-Scatter}
For each $\Phi\in C_c(\mathbf{R}_+^*\times[-1,1])$ and each $h'\in[-1,1]$
\begin{equation}
\label{TransitProba}
\Phi(T_r(h',\cdot))\to\int_0^\infty\int_{-1}^1\Phi(s,h)P(s,h|h')dsdh\quad
	\hbox{ in $L^\infty(\mathbf{S}^1_\omega)$ weak-* as $r\to 0^+$,}
\end{equation}
where the transition probability $P(s,h|h')dsdh$ is the image of the probability 
measure on $K$ given by
\begin{equation}
\label{PartiProba}
d\mu(A,B,Q,N)={\frac{6}{\pi^2}}
\mathbf{1}_{0<A<1}\mathbf{1}_{0<B<1-A}\mathbf{1}_{0<Q<\frac1{2-A-B}}
\frac{dAdBdQ}{1-A}(\delta_{N=0}+\delta_{N=1})
\end{equation}
under the map 
$K\ni (A,B,Q,N)\mapsto\mathbf{T}_{A,B,Q,N}(h')\in\mathbf{R}_+\times[-1,1]$.
Moreover, $P$ satisfies: 
\begin{equation}
\label{PropP}
\begin{array}{r}
\hbox{$(s,h,h')\mapsto(1+s)P(s,h|h')$ is piecewise continuous and bounded
on $\mathbf{R}_+\times[-1,1]\times[-1,1]$,}
\\
\hbox{and $P(s,h|h')=P(s,-h|-h')$Ê for each $h,h'\in[-1,1]$ and $s\ge 0$.}
\end{array}
\end{equation}
\end{theorem}

The proof of (\ref{TransitProba}-\ref{PartiProba}) is based on the explicit 
representation of the transition map in Proposition \ref{P-TransitMap} 
together with Kloosterman sums techniques as in \cite{BocaZaharescu}.
The explicit formula for the transition probability $P$ is very complicated
and we do not report it here, however it clearly entails the properties
(\ref{PropP}).

\section{The Boltzmann-Grad limit of the Lorentz gas dynamics}
\label{S-Dyn}
%%%%%%%%%%%%%%%%%%%%%%%%%%%%%%%%%%%%%%

For each $r\in]0,\frac12[$, denote $d\gamma^+_r(x,\omega)$ the probability 
measure on $\Gamma^+_r$ that is proportional to $\omega\cdot n_xdxd\omega$.
This probability measure is invariant under the billiard map
\begin{equation}
\mathbf{B}_r:\,\Gamma^+_r\ni(x,\omega)\mapsto\mathbf{B}_r(x,\omega)
	=(x+\tau_r(x,\omega)\omega,\mathcal{R}[x+\tau_r(x,\omega)\omega]\omega)
		\in\Gamma^+_r\,.
\end{equation}
For $(x^0,\omega^0)\in\Gamma^+_r$, 
set $(x^n,\omega^n)=\mathbf{B}^n_r(x^0,\omega^0)$ 
and $\alpha^n=\min(|\omega^n_1/\omega^n_2|,|\omega^n_2/\omega^n_1|)$ for each 
$n\ge 0$, and define 
\begin{equation}
b^n_r=(A(\alpha_n,r),B(\alpha_n,r),Q(\alpha_n,r),N(\alpha_n,r)\hbox{ mod. }2)
	\in K
\quad\hbox{  for each $n\ge 0$.}
\end{equation}
We make the following asymptotic independence hypothesis: for each $n\ge 1$ and 
each $\Psi\in C([-1,1]\times K^n)$
$$
\lim_{r\to 0^+}
\int_{\Gamma^+_r}\Psi(h_r,\omega_0,b^1_r,\ldots,b^n_r)d\gamma^+_r(x_0,\omega_0)
=\int_{-1}^1{\textstyle\frac{dh'}2}\int_{\mathbf{S}^1}{\textstyle\frac{d\omega_0}{2\pi}}
\int_{K^n}\Psi(h',\omega_0,\beta_1,\ldots,\beta_n)
	d\mu(\beta_1)\ldots d\mu(\beta_n)\leqno(H)
$$

Under this assumption, the Boltzmann-Grad limit of the Lorentz gas is described by a 
kinetic model on the extended phase space 
$\mathbf{R}^2\times\mathbf{S}^1\times\mathbf{R}_+\times[-1,1]$ --- unlike the Lorentz 
kinetic equation (\ref{LorentzKinEq}), that is set on the usual phase space 
$\mathbf{R}^2\times\mathbf{S^1}$.

\begin{theorem}
\label{T-LimEq}
Assume (H), and let $f^{in}$ be any continuous, compactly supported probability
density on $\mathbf{R}^2\times\mathbf{S}^1$. Denoting by $\tilde R[\theta]$  
the rotation of an angle $\theta$, let $F\equiv F(t,x,\omega,s,h)$ be the solution of
\begin{equation}
\label{EqLim}
\begin{array}{rl}
(\partial_t+\omega\cdot\nabla_x-\partial_s)F(t,x,\omega,s,h)
	&=\int_{-1}^1P(s,h|h')F(t,x,\tilde R[\pi-2\arcsin(h')]\omega,0,h')dh'
\\
F(0,x,\omega,s,h)&=f^{in}(x,\omega)\int_s^\infty\int_{-1}^1P(\tau,h|h')dh'd\tau
\end{array}
\end{equation}
where $(x,\omega,s,h)$ runs through 
$\mathbf{R}^2\times\mathbf{S^1}\times\mathbf{R}^*_+\times]-1,1[$. Then the
family $(f_r)_{0<r<\frac12}$ defined in (\ref{Def-f_r}) satisfies
\begin{equation}
\label{LimBG}
f_r\to\int_0^\infty\int_{-1}^1F(\cdot,\cdot,\cdot,s,h)dsdh
	\hbox{ in $L^\infty(\mathbf{R}_+\times\mathbf{R}^2\times\mathbf{S}^1)$ 
		weak-$*$ as $r\to 0^+$.}
\end{equation}
\end{theorem}

For each $(s_0,h_0)\in\mathbf{R}_+\times[-1,1]$, let $(s_n,h_n)_{n\ge 1}$ be the 
Markov chain defined by the induction formula
\begin{equation}
\label{MarkovChain}
(s_n,h_n)=\mathbf{T}_{\beta_n}(h_{n-1})\hbox{ for each }n\ge 1\,,\quad
\hbox{ together with }\omega_n=\tilde R[2\arcsin(h_{n-1})-\pi]\omega_{n-1}\,,
\end{equation}
where $\beta_n\in K$ are independent random variables distributed under $\mu$. 
The proof of Theorem \ref{T-LimEq} relies upon approximating the particle trajectory
$(X_r,\Omega_r)(t)$ starting from $(x_0,\omega_0)$ in terms of the following jump
process with values in $\mathbf{R}^2\times\mathbf{S^1}\times\mathbf{R}_+\times[-1,1]$
with the help of Proposition \ref{P-TransitMap} 
\begin{equation}
\begin{array}{ll}
(X_t,\Omega_t,S_t,H_t)(x_0,\omega_0,s_0,h_0)=(x_0+t\omega_0,\omega_0,s_0-t,h_0)
	&\qquad\hbox{Êfor }0\le t<s_0\,,
\\
(X_t,\Omega_t,S_t,H_t)(x_0,\omega_0,s_0,h_0)
	=(X_{\tau_n}+(t-s_n)\omega_n,\omega_n,s_{n+1}-t,h_n)
		&\qquad\hbox{for }s_n\le t<s_{n+1}\,.
\end{array}
\end{equation}
Unlike in the case of a random (Poisson) distribution of obstacles, the successive
impact parameters on each particle path are not independent and uniformly distributed
in the periodic case
--- likewise, the successive free path lengths on each particle path are not independent
with exponential distribution. The Markov chain (\ref{MarkovChain}) is introduced to
handle precisely this difficulty.

% etc, etc

% \appendix

% \section{}
% \label{}

% The Acknowledgements are an un-numbered section
%\section*{Acknowledgements}
% Acknowledgements text here

\end{document}